\documentclass[11pt]{article}
\usepackage[DIV=11]{typearea}
\usepackage[T1]{fontenc}
\usepackage{amsmath}
\usepackage{amssymb}
\usepackage{amsthm}

 \DeclareMathOperator{\Ker}{Ker}
 \DeclareMathOperator{\codim}{codim}
 \DeclareMathOperator{\spa}{span}
 \DeclareMathOperator{\sgn}{sgn}
 \newcommand{\eps}{\varepsilon}
 
 \newcommand{\h}{\mathcal{H}}
 \newcommand{\K}{\mathcal{K}}
 \newcommand{\Real}{\mathbb{R}}
 \newcommand{\Complex}{\mathbb{C}}
 \newcommand{\bpola}{|B|^{1/2}}
 \newcommand{\bmin}{|B|^{-1/2}}
 \renewcommand{\l}{\lambda}

 \newcommand{\norm}[1]{\left\Vert#1\right\Vert}
 \newcommand{\knorm}[1]{\norm{#1}_{\K}}
 \newcommand{\deri}{\mathrm{d}}

 \theoremstyle{definition}
 \newtheorem{definition}{Definition}
 \newtheorem{example}{Example}
 \theoremstyle{plain}
 \newtheorem{theorem}{Theorem}
 \newtheorem{lemma}{Lemma} 
 \newtheorem{proposition}{Proposition}
 \newtheorem{corollary}{Corollary}
 \theoremstyle{remark}
 \newtheorem{remark}{Remark}

\begin{document}
\title{On the correspondence between spectra of the operator
pencil $A-\lambda B$ and of the operator $B^{-1}A$}
\author{Ivica Naki\'{c}\thanks{Department of Mathematics, University of Zagreb, 10000 Zagreb, Croatia, e-mail: nakic@math.hr. This work has been partially supported by the HRZZ project grant 9345.
}}

\maketitle

\begin{abstract}
This paper is concerned with the reduction of the spectral problem for symmetric linear operator pencils to a spectral problem for the single operator. 
Also, a Rayleigh--Ritz--like bounds on eigenvalues of linear operator pencils are obtained.


\end{abstract}

\section{Introduction and preliminaries}
\subsection{Introduction}

In many areas of applied mathematics spectral problems for the operator pencil $L(\lambda)=A-\lambda B$ arise, where $A$ and $B$ are symmetric or selfadjoint operators acting in some Hilbert space $\mathcal{H}$. See, for instance, \cite{Dlan}, \cite{Min}, \cite{falaleev2008}, \cite{tretter2000} and references therein. 

One way of dealing with this problem is reducing it to a classical spectral problem for a single operator. 
If $B$ is bounded and boundedly invertible this can be done by introducing the operator $S=B^{-1}A$. 
It is easy to see that the spectral problems for $L$ and $S$ are equivalent.
Additionally, one can think of the operator $S$ as a symmetric or selfadjoint operator in the Krein space $(\mathcal{H},(B\cdot,\cdot))$. 
Here $(\mathcal{H},(B\cdot,\cdot))$ denotes the Banach space $(\mathcal{H},\big\lVert \lvert B \rvert^{1/2} \cdot \big\rVert )$, where $\lvert B \rvert^{1/2} $ is the positive square root of $B$, with the structure of the (possibly indefinite) inner product on this space given by $(B\cdot,\cdot)$.

If $B$ is still bounded, but, say, zero is a point of continuous spectrum of the operator $B$, the space $(\mathcal{H},(B\cdot,\cdot))$ is no longer complete, and $S$ is no longer closed. 
If we drop the assumption that $B$ is bounded, the problems are even more involved, and in the case when zero is an eigenvalue of $B$, $S$ does not exists as an operator.

There exists an extensive literature on the spectral theory of operator pencils with its applications in diverse mathematical and physical domains. Hence, we restrict ourselves to mentioning only two classical works \cite{baumgartel1985} and \cite{marcus1988}. 

The reduction of the spectral problem for operator pencils to the spectral problem for a single operator is a standard procedure used in numerous applications. One of the important steps in the reduction is to show that there is no loss of the relevant spectral information. Sometimes, it is sufficient to show that e.g.\ point spectra coincide, sometimes it is necessary to have a complete correspondence of various types of spectra. There is a plethora of such results in the literature, but usually for a specific setting and a specific type of spectra. Especially important in applications is the point spectrum, and it is treated in several papers, usually in different settings. See, for example \cite{Azizov1990}, \cite{tretter2000}, \cite{lancaster1993}, \cite{cojuhari2007} and \cite{gao2012}. A recent paper \cite{verdier2014} contains some results about the correspondence of spectra in the case of bounded operators between Banach spaces. As far as we are aware, there are no papers which give a systematic treatment of the reduction problem.

This paper treats the case when the operators $A$ and $B$ are symmetric or selfadjoint operators in a Hilbert space $\mathcal{H}$ and the reduction operator is symmetric or selfadjoint in a Krein space $\mathcal{K}$, its topology generated by the operator $B$. 
So the aim is to investigate those reductions for which the symmetry of the problem is preserved. 
We treat various cases depending on the properties of the operator $B$. 

We now give a brief outline of our main results.

In Sections \ref{sec:bounded} and \ref{sec:spectra} we investigate the case when $B$ is bounded and $0\in\sigma_c(B)$. 
In Section \ref{sec:bounded} we give a construction of a reduction operator in the case when $A$ is a symmetric operator with finite and equal deficiency indices.
The main result in Section \ref{sec:bounded} is Theorem \ref{theorem:bijective} where we show that there exists a bijective correspondence between all selfadjoint extensions of the operator $A$ and all selfadjoint extensions of the reduction operator $S$. Section \ref{sec:spectra} is concerned with the correspondence of the spectra. We show that there is a correspondence of spectra and point spectra, but that, in general, there does not exists correspondence for residual and continuous parts of spectra. We give necessary and sufficient conditions for the correspondence of residual, and hence also continuous, parts of spectra. 

In Section \ref{sec:unbounded} we investigate the case when $B$ is unbounded and $0\in\sigma_c(B)$. A construction of the reduction operator follows essentially the one given in \cite{BN}, but the authors there do not investigate the questions we are concerned with, except the case of point spectrum for quasi--uniformly positive operators. The main result in this section is Theorem \ref{theorem:unb} in which we show that under the assumptions stated at the beginning of the section, the point spectrum of the pencil is contained in the point spectrum of the reduction operator and the spectrum of the reduction operator is contained in the spectrum of the pencil. To obtain the correspondence of the spectra, additional assumptions are needed, and some are presented at the end of the section.

In Section \ref{sec:variational}, we apply the results from Sections \ref{sec:bounded} and \ref{sec:unbounded} to obtain a Rayleigh--Ritz--like upper bound for the eigenvalues of the reduction operator given in terms of the eigenvalues of the matrix pair $(A_V,B_V)$, where $A_V$, $B_V$ are orthogonal projections of $A$, $B$ on a finite dimensional subspace. 
To obtain a variational characterization for the reduction operator, we use the results from \cite{BN94}.

In Section \ref{sec:zero} we investigate the case when $0$ is an eigenvalue of $B$. By the use of the theory of linear relations we give constructions of linear relations which act as a substitute for reduced operators. We show that similar results as in Sections \ref{sec:spectra}, \ref{sec:unbounded} and \ref{sec:variational} also hold in this case. Using a result from \cite{S}, we give some examples when a reduced operator can be constructed.

\subsection{Basic definitions and preliminaries}
\label{subsec:basic}

Let $\h$ be a Hilbert space with inner product $(\cdot,\cdot)$ and the corresponding norm $\norm{x}=(x,x)^{1/2}$. 
By $D(T)$ we denote the domain, by $R(T)$ we denote the range of the operator $T$, and by $\Ker(T)$ we denote the null-space of $T$. 
We say that the operator $T$ is symmetric if $T\subset T^{\ast}$ and $T$ has a dense domain. 
The deficiency indices of a symmetric operator $T$ are numbers $n_{\pm}=n_{\pm}(T)$ defined as $n_+=\codim R(T-\mu)$, $n_-=\codim R(T- \overline{\mu})$, where $\mu\in\Complex\setminus\Real$, $\Im\mu>0$ is arbitrary. 
By $\dotplus$ we denote the direct sum.

Let $L(\l)=A-\l B$ be an operator pencil in $\h$, where $A$ and $B$ are densely defined operators in $\h$, such that $B$ is $A$--bounded, i.e.\ $D(B)\supset D(A)$ and $\|Bx\|\le C\|Ax\|$.
Obviously, $D(L(\lambda))=D(A)$.

\begin{definition}
\label{int:spectrum} 
The complex number $\l$ is in the resolvent set of the pencil $L$, (denoted by $\l\in\rho(L)$) if $0\in\rho(L(\l))$. 
The spectrum of $L$ is defined as $\sigma(L)=\Complex \setminus \rho(L)$, i.e.\ $\l\in\sigma(L)$ if $0\in\sigma(L(\l))$. 
Analogously we define the point $\sigma_p$, residual $\sigma_r$, continuous $\sigma_c$, and approximative spectrum $\sigma_{ap}$.
\end{definition}
We also define the set of the \textsl{points of regular type} of the operator pencil $L$ by
\[\Pi(L)=\{\l\in\Complex: \exists C>0 \quad\text{such that}\quad \norm{L(\l)x}\ge C\norm{x}\quad \forall x\in D(L(\l))\},\]
i.e.\ the set of all $\l\in\Complex$ such that $0$ is a point of regular type for the operator $L(\l)$.

We say that the vectors $x_0,x_1,\ldots,x_k\in\h$, $x_0\neq 0$, form a Jordan chain of length $k+1$ for $L$ corresponding to the eigenvalue $\l_0\in\Complex$ if
\[  L(\l_0)x_i=Bx_{i-1},\quad i=0,1,\ldots,k, \]
where $x_{-1}=0$.

In the next proposition we list some properties of linear operator pencils which are essentially known.
\begin{proposition}
\label{int:prop} 
Let $L(\l)=A-\l B$ be an operator pencil in the Hilbert space $\h$, where $B$ is $A$--bounded.
\begin{enumerate}
\item[(i)] If $L$ has a compact resolvent for some $\l\in\rho(L)$, i.e.\ if $L(\l)^{-1}$ is compact, then $L$ has a compact resolvent in each point of the resolvent set.
\item[(ii)] If $B$ is bounded, the set of points of regular type of $L$ is an open set.
\item[(iii)] If $B$ is bounded and $A$ and $B$ are symmetric operators, the deficiency indices of the symmetric operators $A$ and $A-\l B$ coincide, for all $\l\in\Real$.
\end{enumerate}
\end{proposition}
\begin{proof}
The statement (i) easily follows from the generalized resolvent equation. The statement (ii) is evident. For the proof of (iii) see \cite[Theorem 9.3]{GK}.
\end{proof}
The following lemma is known, but we give a proof for reader's convenience.
\begin{lemma}
\label{int:lem} 
Let $\h$ be a Hilbert space, $D$ a dense subspace of $\h$, and $N$ a finite dimensional subspace of $\h$. Then $D\cap N^{\bot}$ is dense in $N^{\bot}$.
\end{lemma}
\begin{proof}\footnote{We thank the anonymous referee for this simplified version of the proof.}
Let $P$ be the orthogonal projection onto $N$. 
Since $P$ is continuous and its range is $N$, the image $P(D)$ is a dense subspace of $N$. As $N$ is finite dimensional, $P(D)=N$. 
Let $D_0$ be a subspace of $D$ such that the restriction $P_0\colon D_0\to N$ of $P$ onto $D_0$ is an isomorphism.
Then $I-P_0^{-1}P$ is continuous on $\mathcal{H}$, its range is $N^{\bot}$ and $(I-P_0^{-1}P)(D)\subseteq D \cap N^{\bot}$.
Since $(I-P_0^{-1}P)(D)$ is dense in $N^{\bot}$, so is $D\cap N^{\bot}$.
\end{proof}

Since we will mainly operate in Krein spaces, we
introduce some basic notions from the Krein space theory. Let $\K$ be a vector space
and let $[\cdot,\cdot]$ be an inner product on $\K$ (not
necessarily definite). The pair $(\K,[\cdot,\cdot])$ is said to be a
Krein space if there exists a direct sum decomposition
$\K=\K_+\dotplus\K_-$ such that $(\K_{\pm},\pm
[\cdot,\cdot])$ are Hilbert spaces and such that $[\K_+,\K_-]=\left\{ 0 \right\}$. For such a
decomposition the corresponding \textsl{fundamental symmetry}
$J$ is a linear operator defined by $J(x_++x_-)=x_+-x_-$,
where $x_+\in\K_+$ and $x_-\in\K_-$. Let
$(\cdot,\cdot)=[J\cdot,\cdot]$. Then the space
$(\K,(\cdot,\cdot))$ is a Hilbert space. Its topology is
independent of the choice of $\K_+$ and $\K_-$. The
subspace $F$ is said to be \textsl{ortho--complemented} if $F\dotplus
F^{\bot}=\K$, where $F^{\bot}=\{x:[x,y]=0\; \forall y\in
F\}$.

The definitions of symmetric and selfadjoint operators in
a Krein space are analogous to those in Hilbert space. A
selfadjoint operator $A$ is said to be \textsl{definitizable} if its
resolvent set is nonempty and there exists a nonzero
polynomial $p$ such that $[p(A)x,x]\ge 0$ for all $x\in
D(p(A))$. A definitizable operator has a spectral function
with finitely many critical points in $\Real
\cup\{\infty\}$. A critical point $z$ is regular if the
spectral function is bounded in a neighborhood of $z$. A
critical point is singular if it is not regular. We say
that $\l$ is an eigenvalue of positive (negative) type of
an operator $T$ if $[x,x]\ge 0$ ($\le 0$) for all eigenvectors
corresponding to $\l$. We say that an operator $T$ is
\textsl{quasi--uniformly positive} (qup) if there exists a
subspace $M$  of finite codimension in $D(T)$ such that
\[   \inf \{[Tx,x]:x\in M,\; \norm{x}=1\}> 0. \]
For more details about Krein spaces see \cite{L},
\cite{Bo} and \cite{AI}. Basic properties of qup operators
can be found in \cite{CN}.

\section{A reduction to the single operator problem for $A$ symmetric and $B$ bounded}
\label{sec:bounded}

In this section we assume:
\begin{enumerate}
\item $A$ is a symmetric operator on a Hilbert space $\mathcal{H}$ with finite and equal deficiency indices $n_+(A)=n_-(A)=:n$,
\item $B$ be is a bounded selfadjoint operator on $\mathcal{H}$ with $0\in\sigma_c(B)$,
\item $\Pi(L)\cap \Real\neq\emptyset$.
\end{enumerate}
Let $E$ be a spectral function of operator $B$. We define
the operator $J$ by $J=E(0,\infty)-E(-\infty,0)$. Operator $J$ is a bounded selfadjoint
operator and $J^2=I$ holds. We define
$|B|:=JB=BJ\ge 0$. Now, we define a space $\hat{\mathcal{K}}$ which consists of sequences $(x_n)$, $x_n\in \mathcal{H}$ in the following way
\[\hat{\mathcal{K}}=\{ (x_n)\colon |B|^{1/2}x_n \text{ converges in } \mathcal{H}\}, \]
i.e.\ as a completion with respect to the norm generated by the operator $\lvert B \rvert^{1/2} $.
We introduce equivalence relation  $\sim$ on $\hat{\mathcal{K}}$:
\[
(x_n)\sim (y_n)\;\text{iff}\;|B|^{1/2}(x_n-y_n)\to 0 \;\text{as}\;
n\to\infty, \]
 and define $\mathcal{K}$ as the corresponding quotient space
$\mathcal{K}=\hat{\mathcal{K}}/{\sim}$. On $\mathcal{K}$ we
introduce two inner products $(\cdot,\cdot)_{\K}$ and
$[\cdot,\cdot]_{\K}$ as follows:
\begin{gather*}
((x_n),(y_n))_{\K}=\lim_{n\to\infty}(|B|^{1/2}x_n,|B|^{1/2}y_n),
\\
[(x_n),(y_n)]_{\K}=\lim_{n\to\infty}(J|B|^{1/2}x_n,|B|^{1/2}y_n).
\end{gather*}
It is easy to see that $(\K,(\cdot,\cdot)_{\K})$ is a Hilbert
space, and that $(\K,[\cdot,\cdot]_{\K})$ is a Krein space. These
spaces have the norm given by $\norm{(x_n)}_{\K}=\lim_{n\to\infty}
\norm{|B|^{1/2}x_n}$. If we identify $x\in\h$ with the sequence
$(x,x,\ldots)$ in $\hat{\mathcal{K}}$, it is evident that $\h$ can
be regarded as a subspace in $\K$, and in that case
$\knorm{x}=\norm{\bpola x}$ holds. Also, $\h$ is a dense subspace in
$\K$. Operators in $\h$ can also be regarded as operators in
$\K$. Also, $J$ can be extended to a bounded operator in $\K$, and
this extension we also denote with $J$. It is easy to see
that $J$ is a fundamental symmetry in $\K$.
\begin{proposition}
\label{prop:isom} The spaces $\h$ and $\K$ are isomorphic.
\end{proposition}
\begin{proof}
We define
\begin{equation}
 	\label{eq:defT}
 	T:\K\to\h,\quad T(x_n)=\lim_{n\to \infty}\bpola x_n.
 \end{equation} 
Then it is easy to see that $T$ is a linear surjection from $(\K,(\cdot,\cdot)_{\K})$ to $(\h,(\cdot,\cdot))$ which preserves the inner products, i.e.\ $(T(x_n),T(y_n))=((x_n),(y_n))_{\K}$.
\end{proof}
In applications, the space $\K$ can usually be represented
as a function space as can be seen from the following example.
\begin{example}
Let $\h=L_2(a,b)$, and let $w\in L_{\infty}(a,b)$ be a function which is zero on some non--empty set of measure zero.
We define the operator $B$ by $(Bf)(x)=w(x)f(x)$, $f\in L_2(a,b)$. 
The operator $B$ is obviously bounded and $0\in \sigma_c(B)$ holds.

The Krein space $\K$ can be identified with the space
\[\{ f: \sqrt{|w(\cdot)|}f\in L_2(a,b)\}\] (i.e.\ a weighted $L_2$ space with the weight $w$, see \cite[Section 8.4]{Wei}) with the inner products given by
\begin{gather*}
[f,g]_{\K}=\int_{a}^b w(x) f(x)\overline{g(x)}\deri x,\\
(f,g)_{\K}=\int_{a}^b |w(x)| f(x)\overline{g(x)}\deri x.
\end{gather*}
The operator $T$ is given by $Tf=\sqrt{|w(\cdot)|}f$, and the fundamental symmetry $J$ is given by $Jf = (\sgn w)f$.
\end{example}
\begin{lemma}
\label{lem:dense} Let $D$ be a dense subset in $\h$. Then
$D$ is also dense in $\K$.
\end{lemma}
\begin{proof}
Let $\eps >0, x\in\h$ be arbitrary. Since $R(\bpola)$ is
dense in $\h$, there exists $y\in\h$ such that
$\norm{\bpola y-x}<\frac{\eps}{2}$, and since $D$ is dense
in $\h$, there exists $z\in D$ such that
$\norm{y-z}<\frac{\eps}{2}\norm{\bpola}^{-1}$. Now
\begin{equation*}
\norm{\bpola z-x}\le\norm{ \bpola (z-y)}+\norm{\bpola
y-x}<\eps.
\end{equation*}
This implies that $\{\bpola x: x\in D\}$ is dense in $\h$.
Now let $(x_n)\in\K, \\
 \eps >0$ be arbitrary. Set $\widetilde x=T(x_n)\in\h$. Then there
is $y\in\h$ such that $\norm{\widetilde x -\bpola y}<\eps$,
which implies $\knorm{(x_n) -y}<\eps$.
\end{proof}
\begin{proposition}
\label{prop:symm} The operator $B^{-1}A$ is symmetric in
$(\K,[\cdot,\cdot]_{\K})$.
\end{proposition}
\begin{proof}
Let $x,y\in D(B^{-1}A)=\{x\in D(A):Ax\in R(B)\}$. Then
\begin{equation*}
\begin{split}
[B^{-1}Ax,y]_{\K}& =(J\bpola B^{-1}A x,\bpola y)=(Ax,y) \\
&= (x,Ay)=(J\bpola x,\bpola B^{-1}Ay)=[x,B^{-1}Ay]_{\K}.
\end{split}
\end{equation*}
Let $\l\in\Pi(L)\cap\Real$ be arbitrary. Then $R(L(\l))$ is closed.
From Proposition \ref{int:prop} and
\cite[Theorem 9.2]{GK} follows $\codim R(L(\l))=n$.
 From Lemma \ref{int:lem}
follows that $R(L(\l))\cap R(B)$ is dense in $R(L(\l))$.
Note that
\begin{equation*}
\begin{split}
 R(L(\l))\cap R(B)&=\{ (A-\l B)x: x\in D(A), Ax\in R(B)\} \\
 &=\{(A-\l B)x: x\in
D(B^{-1}A)\}.
\end{split}
\end{equation*}
Let $C>0$ be such that $\norm{L(\l)x}\ge C\norm{x}$. Let
$x\in D(A)$, $\eps >0$ be arbitrary. Set $\widetilde
x=L(\l)x$. Then it exists $y\in D(B^{-1}A)$ such that
$\norm{L(\l)y-\widetilde x}<C\eps$. Then
\begin{equation*}
C\eps > \norm{L(\l)y-\widetilde x}=\norm{L(\l)(y-x)}\ge
C\norm{x-y},
\end{equation*}
hence $\norm{x-y}<\eps$. So $D(B^{-1}A)$ is dense in
$D(A)$, which implies that $D(B^{-1}A)$ is dense in $\h$.
Now Lemma \ref{lem:dense} implies that
$D(B^{-1}A)$ is dense in $\K$, hence $B^{-1}A$ is a
symmetric operator in $\K$.
\end{proof}
We will need the following three lemmas.
\begin{lemma}
\label{lem:xunulu} Let $(x_n)\in D(A^{\ast})$ be a
sequence such that
\begin{equation}
\label{eq:ast} A^{\ast}x_n\to 0
\end{equation}
and
\begin{equation}
\label{eq:bpol} \bpola x_n\to 0.
\end{equation}
Then $x_n\to 0$.
\end{lemma}
\begin{proof}
From (\ref{eq:bpol}) it follows that $Bx_n\to 0$. Let $\lambda\in\Pi(L)\cap \Real$ be arbitrary.
Since $(A-\lambda B)^{\ast}=A^{\ast}-\lambda B$, we have
\[ D(A^{\ast})=D(A)\dotplus\Ker (A^{\ast}-\lambda B). \]
Let $x_n=x_n^0+y_n$, $x_n^0\in D(A)$, $y_n\in\Ker (A^{\ast}-\lambda B)$.
Since $(A^{\ast}-\lambda B)x_n=(A-\lambda B)x_n^0$, and $(A-\lambda B)^{-1}$ is a bounded operator, it follows $x_n^0\to 0$. 
Hence $Bx_n^0\to 0$ and $By_n= B(x_n -x_n^0)\to 0$.
Since the operator
\[B|_{\Ker(A^{\ast}-\lambda B)}:\Ker(A^{\ast}-\lambda B)\to B(\Ker(A^{\ast}-\lambda B))\]
acts between finite--dimensional spaces and is non--singular, we have $y_n\to 0$, hence $x_n\to 0$.
\end{proof}
\begin{lemma}
\label{lem:abc}
We have
\begin{enumerate}
\item[(i)] The set $D(\bmin A^{\ast })$ is a core of $A^{\ast}$.
\item[(ii)] The set $D(\bmin A)$ is a core of $A$.
\item[(iii)] For each selfadjoint extension $\widetilde A$ of $A$, the set $D(\bmin \widetilde{A})$ is a core of $\widetilde A$.
\end{enumerate}
\end{lemma}
\begin{proof}
Let $\l\in\Pi(L)\cap \Real$ be arbitrary.

(i) As in the proof of Proposition
\ref{prop:symm}, we have $\codim R(A-\l B)=n$. Also
$(A-\l B)^{\ast}=A^{\ast}-\l B$ holds. This implies
$D(A^{\ast})=D(A)\dotplus\Ker (A^{\ast}-\l B)$.

Let $x\in
D(A^{\ast})$ be arbitrary, and let $x=x_0+y_0$, $x_0\in
D(A)$, $y_0\in\Ker(A^{\ast}-\l B)$. Since $\Ker(A^{\ast}-\l
B)\subset D(\bmin A^{\ast})$, it is enough to find a
sequence $(x_n^0)\in D(|B|^{-1/2}A)$ such that $x_n^0\to x_0$ and $Ax_n^0\to
Ax_0$. Lemma \ref{int:lem} implies that $R(A-\l B)\cap
R(\bpola)$ is dense in $R(A-\l B)$, so for $(A-\l B)x_0$
there exists a sequence $x_n^0\in D(\bmin A)$ such that
$(A-\l B)x_n^0\to (A-\l B)x_0$. Since $(A-\l B)^{-1}$ is
bounded it follows that $x_n^0\to x_0$  which implies
$Bx_n^0\to Bx_0$, and consequently, $Ax_n^0\to Ax_0$.

(ii) Since $R(A-\lambda B)\cap R(B)$ is dense in $R(A-\lambda B)$, for an arbitrary $x\in D(A)$ there exists a
sequence $x_n\in D(B^{-1}A)$ such that $(A-\lambda B)x_n\to (A-\lambda B)x$, which implies $x_n\to x$ and $Ax_n\to Ax$.

(iii) Since $ \widetilde{A}$ is a finite--dimensional extension of $A$, $ \mathrm{codim}R( \widetilde{A}-\lambda B)\le n$, and $R( \widetilde{A}-\lambda B)$ is closed. 
Hence, again from Lemma \ref{int:lem}, $R( \widetilde{A}-\lambda B)\cap R(|B|^{1/2})$ is dense in $R( \widetilde{A}-\lambda B)$. 
Now we can proceed as in the proof of the first statement of the Lemma.
\end{proof}
\begin{lemma}
\label{lem:adj} The formulae $(\bmin A)^{\ast}=A^{\ast}\bmin$ and $(B^{-1} A)^{\ast}=A^{\ast}B^{-1}$ hold.
\end{lemma}
\begin{proof}
Since $\bpola$ is bounded, \cite[Satz 4.19]{Wei} implies
\begin{equation}
\label{eq:AA}
(\bpola (\bmin A))^{\ast}=(\bmin A)^{\ast}\bpola.
\end{equation}
From Lemma \ref{lem:xunulu} follows
\begin{equation}
\label{eq:BB}
(\bpola (\bmin A))^{\ast}=(A|_{D(\bmin A)})^{\ast}
=(\overline{A|_{D(\bmin A)}})^{\ast}=A^{\ast},
\end{equation}
Now (\ref{eq:AA}) and (\ref{eq:BB}) imply $(\bmin A)^{\ast}\bpola=A^{\ast}$.
The second formula follows analogously.
\end{proof}
Let us define the operator $S$ in $\K$ by
\begin{gather*}
 D(S)=\{x\in D(A):Ax\in R(\bpola)\}=
 D(\bmin A) \\
 Sx=T^{-1}J\bmin Ax,
\end{gather*}
where the operator $T$ is defined by \eqref{eq:defT}.
\begin{proposition}
\label{prop:Sclosed}
The operator $S$ is closed.
\end{proposition}
\begin{proof}
Let $(x_n)$ be a sequence in $D(S)$, and let
$x_n \to x$ in the norm $\knorm{\cdot}$, which implies $x\sim
(x_n)$. Set
\begin{equation}
\label{eq:Tx} \widetilde x=T(x_n),
\end{equation}
and let $Sx_n\to y$ in the norm of $\K$. It follows
\begin{equation}
\label{eq:ytilda} J\bmin A x_n\to \widetilde y:=Ty.
\end{equation}
From (\ref{eq:Tx}) and (\ref{eq:ytilda}) follows
$Bx_n\to J\bpola \widetilde x$ and $Ax_n\to J\bpola
\widetilde y$, \\ respectively, which implies
\begin{equation}
\label{eq:a-lb} (A-\l B)x_n\to J\bpola (\widetilde y -\l
\widetilde x ),
\end{equation}
for all $\l\in\Complex$. Especially, (\ref{eq:a-lb})
holds for $\l_0\in\Pi(L)$, which implies
\[x_n\to (A-\l_0
B)^{-1}J\bpola (\widetilde y-\l_0 \widetilde x)=:x_0. \]
Hence we found $x_0\in \h$ such that
$x_n \to x_0$ in the norm $\norm{\cdot}$ and $Ax_n\to
J\bpola \widetilde y$, which implies $x_0\in D(A)$ and
$Ax_0=J\bpola \widetilde y$. Hence $Ax_0\in R(\bpola)$,
and $J\bmin A x_0=\widetilde y$, i.e.\ $Sx_0=y$ and $x_0\in
D(S)$.
\end{proof}
\begin{proposition}
The adjoint of $S$ in $\K$ is given by
\begin{gather*}
 D(S^{\ast})=\{x\in D(A^{\ast}):A^{\ast}x\in R(\bpola)\}=
 D(\bmin A^{\ast}) \\
 S^{\ast}x=(y_n), \quad\text{where}\;\; T(y_n)=J\bmin A^{\ast} x.
\end{gather*}
\end{proposition}
\begin{proof}
Let $(y_n),(z_n)\in\K$ be such that
\[
[Sx,(y_n)]_{\K}=[x,(z_n)]_{\K},\quad\forall x\in D(S).
\]
This can be written as
\begin{equation}
\label{eq:adj}
 \lim_n (\bmin Ax,\bpola y_n)=\lim_n (J\bpola x,\bpola z_n)
 \quad\forall x\in D(S).
\end{equation}
Let us denote $\widetilde{y}=T(y_n)$,
$\widetilde{z}=T(z_n)$. Then (\ref{eq:adj}) reads
\[ (\bmin Ax, \widetilde{y})=(x,J\bpola \widetilde{z})\quad
\forall x\in D(\bmin A).\]
From Lemma \ref{lem:adj} follows $\widetilde{y}\in
D(A^{\ast}\bmin)$ and
\[   A^{\ast}\bmin \widetilde{y}=J\bpola \widetilde{z}. \]
Now, since $\widetilde{y}\in R(\bpola)$ there exists
$y\in\h$ such that $(y_n)\sim y$. Hence we can take
$\widetilde{y}=\bpola y$, which implies $A^{\ast}y=J\bmin
\widetilde{z}$. Hence $J\bmin A^{\ast}y=\widetilde{z}$, which implies
$y\in D(\bmin A^{\ast})$ and $T(z_n)=J\bmin A^{\ast}y$.

On the other hand, take an arbitrary $y\in  D(\bmin A^{\ast})$. Then one can easily prove
\[ (\bmin Ax,\bpola y)=(\bpola x,\bmin A^{\ast} y), \; \forall x\in D(S),\]
hence $y\in D(S^{\ast})$ and $T(S^{\ast}y)=J\bmin A^{\ast}y$.
\end{proof}
\begin{remark}
If $n=0$, i.e.\ if $A$ is selfadjoint, then $S$ is also
selfadjoint.
\end{remark}
\begin{theorem}
\label{theorem:S-op}
The operator $S$ is the closure of the operator $B^{-1}A$ in $\K$.
\end{theorem}
\begin{proof}
Let us denote with $S'$ the closure of the operator $B^{-1}A$ in $\K$.

Let $x\in D(B^{-1}A)$, and set $(y_n)=Sx$, i.e.\ $\lim_{n\to \infty}
\bpola y_n=J\bmin Ax\in R(\bpola)$. It follows that
there exists $y\in \h$ such that $y\sim (y_n)$, and
$\bpola y=J\bmin Ax $, hence $y=B^{-1}Ax$. Hence we have proved $B^{-1}A
\subset S$, so $S'\subset S$.

Let $(y_n)\in D((B^{-1}A)^{\ast})$ be arbitrary, where the adjoint is taken with respect to the
geometry in $\K$, and set $(z_n)=(B^{-1}A)^{\ast}x$.
Then
\[ [B^{-1}Ax,(y_n)]_{\K}=[x,(z_n)]_{\K} \]
holds for all $x\in D(B^{-1}A)$.
This implies
\begin{equation}
\label{eq:Sstar}
(\bmin Ax,\widetilde{y})=(x,J\bpola \widetilde{z}),\; \forall x\in D(B^{-1}A),
\end{equation}
where $\widetilde{y}=T(y_n)$, $ \widetilde{z}=T(z_n)$.
The relation (\ref{eq:Sstar}) can be written as
\begin{equation}
\label{eq:CC}
(B^{-1}Ax,J\bpola \widetilde{y})=(x,J\bpola \widetilde{z}).
\end{equation}
From Lemma \ref{lem:adj} we know that $(B^{-1}A)^{\ast}=A^{\ast}B^{-1}$, hence (\ref{eq:CC}) implies
\[ A^{\ast}\bmin \widetilde{y}=J\bpola \widetilde{z}. \]
This implies $ \widetilde{y}\in R(\bpola)$, hence $(y_n)\sim y\in \h$ and $A^{\ast}y=J\bpola \widetilde{z}$,
which implies $y\in D(\bmin A^{\ast})$. We have proved that $(B^{-1}A)^{\ast}\subset S^{\ast}$, hence
$S\subset (B^{-1}A)^{\ast\ast}=S'$.
\end{proof}
Let $S$ be a closed symmetric operator in a Hilbert space with equal deficiency indices. 
It is well-known that the set $D(S)$ provided with the inner product $(x,y)_{D(S^{\ast})}=(x,y)+(S^{\ast}x, S^{\ast}y)$ is a Hilbert space.
\begin{definition}[\cite{schmudgen2012},\cite{DM},\cite{GG}]
Let $S$ be a symmetric operator in a Hilbert space $\h$. 
A collection $\{\Omega,\Gamma_1,\Gamma_2\}$, in which $\Omega$ is a Hilbert space, and $\Gamma_1, \Gamma_2:D(C^{\ast})\to\Omega$ are bounded operators, is called a space of boundary values (SBV) for $S^{\ast}$, if
\begin{enumerate}
\item $(S^{\ast}x,y)-(x,S^{\ast}y)=(\Gamma_1x,\Gamma_2y)_{\Omega}-
(\Gamma_2x,\Gamma_1y)_{\Omega}\quad \forall x,y\in
D(S^{\ast}),$
\item the mapping $\Gamma:x\mapsto\{\Gamma_1x,\Gamma_2x\}$
from $D(S^{\ast})$ to
$\Omega\times \Omega$ is surjective.
\end{enumerate}
\end{definition}
In the case of a symmetric operator $S$ in a Krein space we have
\begin{definition}[\cite{Ko}]
\label{def:ksbv} Let $S$ be a symmetric operator in a Krein space $(\K,[\cdot,\cdot])$. 
A collection $\{\Omega,\widehat{J},\Gamma_1,\Gamma_2\}$, in which $\Omega$ is a Hilbert space, $\widehat{J}$ is an
involution in $\Omega$, and $\Gamma_1,\Gamma_2: D(S^{\ast})\to \Omega$ are bounded operators, is called a space of boundary values for $S^{\ast}$, if
\begin{enumerate}
\item $ [S^{\ast}x,y]-[x,S^{\ast}y]=(\Gamma_1x,
\widehat{J}\Gamma_2y)_{\Omega}-
(\Gamma_2x, \widehat{J}\Gamma_1y)_{\Omega}\quad\forall
x,y\in D(S^{\ast}),$
\item the mapping $\Gamma:x\mapsto\{\Gamma_1x,\Gamma_2x\}$
from $D(S^{\ast})$ to
$\Omega\times \Omega$ is surjective.
\end{enumerate}
\end{definition}
In the case of Hilbert space operators, it can be shown
that a SBV exists for any symmetric operator with equal
deficiency indices (see \cite{DM}, \cite{GG}). In the case
of Krein space operators, it can be shown that a SBV
exists for an operator $T$ if $\Pi(T)\neq\emptyset$ holds,
and if the deficiency indices of $T$ are equal (see \cite{Ko}).
For a fixed SBV for the operator $T$, there exists a
bijective correspondence between the collection of the
closed extensions $\widetilde{T}$ and the set of closed
relations $\theta$ in $\Omega$:
\begin{equation}
\label{eq:ZZ}
x\in D(\widetilde{T})=D(\widetilde{T}_{\theta})
\Longleftrightarrow
\{\Gamma_2x,\Gamma_1x\}\in\theta\subset\Omega\times\Omega.
\end{equation}
An extension $\widetilde{T}$ is selfadjoint  if and
only if $\theta$ is a selfadjoint relation, or
equivalently, if there exists a selfadjoint operator $G$ on
$\Omega$ given by the relation
\[ (\cos G)\Gamma_2x-(\sin G)\Gamma_1x=0 \Longleftrightarrow \{\Gamma_2x,\Gamma_1x\} \in\theta. \] 
For more details on the theory of SBV see \cite{GG}, \cite{Ko}, \cite{DM} and \cite{schmudgen2012}.

Let $\{\Omega,\Gamma_1,\Gamma_2\}$ be a fixed SBV for the operator
$A^{\ast}$. Set
$\widetilde{\Gamma}_i=\Gamma_i|_{D(\bmin A^{\ast})}$,
$i=1,2$. We will show that a collection $\{\Omega, I,
\widetilde{\Gamma}_1,\widetilde{\Gamma}_2\}$ is a SBV for
the operator $S^{\ast}$. Hence, we need to show that
the statements 1. and 2. from Definition \ref{def:ksbv} hold. We start with
\begin{equation*}
\begin{split}
[S^{\ast}x,y]_{\K}-[x,S^{\ast}y]_{\K} &= (\bmin
A^{\ast}x,\bpola y)- (J\bpola x, J \bmin A^{\ast}y) \\
&=(A^{\ast}x,y)-(x,A^{\ast}y)=(\Gamma_1x,\Gamma_2y)_{\Omega}-
(\Gamma_2x,\Gamma_1y)_
{\Omega} \\
&=(\widetilde{\Gamma}_1x,\widetilde{\Gamma}_2y)_{\Omega}-
(\widetilde{\Gamma}_2x,\widetilde{\Gamma}_1y)_{\Omega}.
\end{split}
\end{equation*}
To see that $\widetilde{\Gamma_i}$ is bounded for $i=1,2$,
it is sufficient to prove that there exists $C>0$ such that
\[ \|x\|_{D(A^{\ast})}\le C \|x\|_{D(S^{\ast})}\quad \forall
x\in D(S^{\ast}), \]
and since
\[ \norm{A^{\ast}x}\le\norm{\bpola \bmin A^{\ast}x}\le
\norm{\bpola}\norm{\bmin
A^{\ast}x}, \] it is enough to see that there exists $C>0$
such that
\begin{equation}
\label{eq:cba} \norm{x}\le C\left(\norm{\bpola
x}+\norm{A^{\ast}x}\right) \quad\forall x\in D(\bmin A^{\ast}).
\end{equation}
But this is implied by Lemma \ref{lem:xunulu}. Indeed,
suppose that (\ref{eq:cba}) does not hold. Then there is a
sequence $(x_n)\in D(\bmin A^{\ast})$, $\norm{x_n} =1$,
such that
\[  \norm{x_n}\ge n\left(\norm{\bpola x_n}+\norm{A^{\ast}x_n}\right), \]
which implies $A^{\ast}x_n\to 0$ and $\bpola x_n \to 0$,
hence $x_n\to 0$, which is a contradiction with the
statement of Lemma \ref{lem:xunulu}.

The density of $D(S^{\ast})$ in $D(A^{\ast})$ follows
directly from Lemma \ref{lem:abc}. Our next aim is
to show that $\widetilde{\Gamma}=
\{\widetilde{\Gamma}_1,\widetilde{\Gamma}_2\}$ is
surjective. Let $f,g\in\Omega$ be arbitrary. Then there
exists $x\in D(A^{\ast})$ such that $\Gamma_1 x=f$,
$\Gamma_2 x=g$. Since $D(S^{\ast})$ is dense in
$D(A^{\ast})$, there exists a sequence $(x_n)\in
D(S^{\ast})$ such that $x_n\to x$ in the norm of
$D(A^{\ast})$, and since $\Gamma_i$, $i=1,2$ are bounded, it
follows that $\widetilde{\Gamma}_ix_n= \Gamma_i x_n\to \Gamma_i
x$, $i=1,2$, so $R(\widetilde{\Gamma})$ is dense in
$\Omega\times\Omega$. Now, since the space
$\Omega\times\Omega$ is finite--dimensional,
 it follows $R(\widetilde{\Gamma})=\Omega\times
\Omega$.

Hence, we can parameterize the selfadjoint extensions of the operator $S$ by the use of a SBV for
the operator $A^{\ast}$.
More precisely, let $\widetilde{A}=\widetilde{A}_{\theta}$ be an extension
of the operator $A$ generated by the closed relation $\theta$ on
$\Omega$. We define an extension $\widetilde{S}$ of $S$ in
the following way:
\[ D(\widetilde{S})=D(\bmin \widetilde{A})\quad\text{and}\quad
\widetilde{S}x=S^{\ast}x. \] This definition implies
\begin{equation*}
\begin{split}
D(\widetilde{S})&=\{x\in D(\widetilde{A}):
\widetilde{A}x=A^{\ast}x\in R(\bpola)\} \\
&=\left\{x\in D(A^{\ast}):A^{\ast}x\in
R(\bpola):\{\Gamma_2x,\Gamma_1x\}\in\theta\right\} \\
&=\left\{x\in D(S^{\ast}):\{\widetilde{\Gamma}_2x,\widetilde{\Gamma}_1x\}\in
D(\widetilde{S}_{\theta})\right\}.
\end{split}
\end{equation*}
On the other hand, if
$\widetilde{S}=\widetilde{S}_{\theta}$ is a closed
extension of the operator $S$ generated by a closed relation
$\theta$ on $\Omega$, then it generates a closed extension
$\widetilde{A}_{\theta}$ of the operator $A$. Thus, we have shown
\begin{theorem}
\label{theorem:bijective} There exists a bijective
correspondence between all closed extensions of the
operator $A$ and all closed extensions of the operator $S$.
Specifically, there exists a bijective correspondence
between all selfadjoint extensions of the operator $A$ and
all selfadjoint extensions of the operator $S$. Moreover,
if $\widetilde{A}$ is a selfadjoint extension of $A$, then
the operator $\widetilde{S}$ defined by
$D(\widetilde{S})=D(\bmin \widetilde{A})$, $\widetilde{S}x=(y_n)$,
$T(y_n)=J\bmin \widetilde{A}$, is a selfadjoint extension
of $S$.
\end{theorem}
\begin{remark}
In the case when a selfadjoint extension of $A$ is given by the
relation $(\cos G)\Gamma_2x-(\sin G)\Gamma_1x=0$, where
$G$ is a selfadjoint operator in $\Omega$, the relation
$(\cos G)\widetilde{\Gamma}_2x-(\sin G)\widetilde{\Gamma}_1
x=0$ defines a selfadjoint extension of $S$.
\end{remark}
\begin{example}
Let $\h=L_2(0,1)$, $Af=-\frac{\deri^2f}{\deri x^2}$, $D(A)=C^2(0,1)$, $(Bf)(x)=xf(x)$.
The operators $A$ and $B$ are obviously symmetric, $B$ is bounded and $0\in \sigma_c(B)$.
One can easily see that $\Pi(L)\ne\emptyset$.
A SBV for $A^{\ast}$ can be chosen as follows:
\[ \Omega=\Complex^2,\; \Gamma_1=\{-f(0),f(1)\},\; \Gamma_2=\{f'(0),f'(1)\}. \]
We have $D(A^{\ast})=W_2^2(0,1)$, where $W_2^2(0,1)$ is the usual Sobolev space, and
\[ D(\bmin A^{\ast})=\left\{ f\in W_2^2(01,): \frac{1}{\sqrt{x}}\frac{\deri^2f}{\deri x^2}(x)\in L_2(0,1)\right\}. \]
The collection $\{\Complex^2,I,{\Gamma_1|}_{D(\bmin A^{\ast})},{\Gamma_2|}_{D(\bmin A^{\ast})}\}$ is a SBV for $S^{\ast}$,
hence all selfadjoint extensions of $S$ are parameterized by the selfadjoint relations in $\Complex^2$ by the formula
(\ref{eq:ZZ}).
\end{example}
\section{The correspondence of the spectra in the case of bounded $B$}
\label{sec:spectra}

In this section we have the same assumptions as in the previous section. 
Our aim in this section is to show the connection between the spectra of $\widetilde{L}$ and $\widetilde{S}$.

Let $\theta$ be a selfadjoint relation in
$\Omega\times\Omega$ such that
$\rho(\widetilde{L})\neq\emptyset$ holds, where
$\widetilde{L}(\l)=\widetilde{A}_{\theta}-\l
B=\widetilde{A}-\lambda B$, and let
$\widetilde{S}=\widetilde{S}_{\theta}$ be a selfadjoint
extension generated by the relation $\theta$.
\begin{proposition}
The operator $\widetilde{S}$ is not bounded.
\end{proposition}
\begin{proof}
Let us assume that $\widetilde{S}$ is bounded. This
implies $D(\widetilde{S})= \K$, hence $D(\bmin
\widetilde{A})=\K$, which implies $\K=\h$. This implies that
for each sequence $(x_n)\in\h$ such that $\bpola x_n\to
y$, there exists $x\in\h$ such that $x\sim (x_n)$, i.e.\
$\bpola x=y$. Hence $R(\bpola)=\h$, a contradiction with
the assumption $0\in\sigma_c(B)$.
\end{proof}
For a symmetric operator $S$ in a Hilbert or Krein space,
by $\pi(S)$ we denote the number of negative squares of
Hermitian sesquilinear form $(Sx,y)$, $x,y\in D(s)$, where
$(\cdot,\cdot)$ is the inner product on the Hilbert or Krein
space. 
That is, $\pi(S)$ is the supremum of the dimensions of all subspaces $\mathcal{L}$ such that $(Sx,x)<0$ for all $0\ne x \in \mathcal{L}$.
\begin{theorem}
\label{theorem:kva}
For each selfadjoint extension $\widetilde{A}$ of $A$, and
$\widetilde{S}$ the corresponding extension of $S$, the formula
$\pi(\widetilde{S})=\pi(\widetilde{A})$
holds.
\end{theorem}
\begin{proof}
The inequality $\pi(\widetilde{S})\le\pi(\widetilde{A})$
is evident.

We will first prove the theorem in the case
$\pi(\widetilde{S})=0$. In this case we have
$[\widetilde{S}x,x]_{\K}\ge 0$ for all $x\in
D(\widetilde{S})$, or, equivalently
$(\widetilde{A}x,x)\ge0$ for all $x\in D(\bmin
\widetilde{A})$. Let us assume that
$\pi(\widetilde{A})>0$. Then there exists some $x\in
D(\widetilde{A})$ such that $(\widetilde{A}x,x)<0$.
Lemma \ref{lem:abc} implies that there is a
sequence $(x_n)\in D(\bmin \widetilde{A})$ such that
$x_n\to x$ and $\widetilde{A}x_n\to \widetilde{A}x$, hence
$(\widetilde{A}x_n,x_n)\to (\widetilde{A}x,x)$, so for $n$
large enough we have $(\widetilde{A}x_n,x_n)<0$; a
contradiction.

Now we treat the case $\pi(\widetilde{S})>0$. Let integer
$n$ and elements $x_1,\ldots,x_n\in D(\widetilde{A})$ be
arbitrarily chosen. Let $\mathcal{A}$ be the matrix
$((\widetilde{A}x_n,x_n))_{i,j=1}^n$, and let $k$ be the
number of negative eigenvalues of $\mathcal{A}$. Then, by Lemma \ref{lem:abc} (iii), for
each $\eps>0$ there exist $x_1',\ldots,x_n'\in D(\bmin
\widetilde{A})$ such that $\norm{x_i-x_i'}<\eps $ and
$\|\widetilde{A}(x_i-x_i')\|<\eps$. Since
\begin{equation*}
\begin{split}
|(\widetilde{A}x_i,x_j)-(\widetilde{A}x_i',x_j')|&=
|(\widetilde{A}(x_i-x_i'),x_j)+(\widetilde{A}x_i',x_j-x_j')|<
\\ &< \eps (\norm{x_j}+\norm{Ax_i})+\eps^2,
\end{split}
\end{equation*}
for a sufficiently small $\eps >0$ the matrix
$\mathcal{A}'= ((\widetilde{A}x_i',x_j'))_{i,j=1}^n$ will
have at least $k$ negative eigenvalues. Hence, the number
of negative squares of the Hermitian form $(\widetilde{A}x,y)$,
$x,y\in D(\widetilde{A})$ is smaller or equals the number
of negative squares of the Hermitian form
$[\widetilde{S}x,y]_{\K}$, $x,y\in D(\widetilde{S})$, which
finishes the proof.
\end{proof}
\begin{theorem}
\label{theorem:spect}
The spectra of $\widetilde{S}$ and
$\widetilde{L}$ coincide. Moreover, their point spectra
coincide, as well as their eigenspaces and Jordan chains.
\end{theorem}
\begin{proof}
First we will show
$\sigma_p(\widetilde{L})=\sigma_p(\widetilde{S})$. Let
$\l\in\sigma_p(\widetilde{L})$ be arbitrary. Then there exists $x\in
D(\widetilde{A})$, $x\neq 0$ such that
\begin{equation}
\label{eq:atilda=lb} \widetilde{A}x=\l Bx.
\end{equation}
It is evident that $x\in D(\widetilde{S})$. Multiplying
(\ref{eq:atilda=lb}) by $J\bmin$ we obtain \\
$J\bmin\widetilde{A}x=\l\bpola$, which implies
$T(\widetilde{S}-\l)x=0$, hence $(\widetilde{S}-\l)x=0$.

On the other hand, for $\l\in\sigma_p(\widetilde{S})$
there exists $x\in D(\widetilde{S})$, $x\neq 0$, such that
\begin{equation}
\label{eq:ps} J\bmin \widetilde{A}x=\l\bpola x.
\end{equation}
Multiplying (\ref{eq:ps}) by $J\bpola$ we obtain
$\widetilde{A}x=\l Bx$, i.e.\
$\l\in\sigma_p(\widetilde{L})$. Hence,
$\sigma_p(\widetilde{L})= \sigma_p(\widetilde{S})$.

The relation $(\widetilde{S}-\l)x_1=x_0$ implies
\begin{equation}
\label{eq:jl} J\bmin \widetilde{A}x_1-\l \bpola x_1=\bpola
x_0.
\end{equation}
Multiplying (\ref{eq:jl}) by $J\bpola$ we obtain
$(\widetilde{A}-\l B)x_1=Bx_0$, hence a Jordan chain of
$\widetilde{S}$ corresponding to an eigenvalue $\l$ is
also a Jordan chain of $\widetilde{L}$ for the same
eigenvalue $\l$. The other direction can be seen
analogously.

Let $\l\in\rho(\widetilde{S})$. This implies
$T(R(\widetilde{S}-\l))=\h$, i.e.\\ $R(J\bmin
\widetilde{A}-\l\bpola)=\h$. This implies
\begin{equation}
\label{eq:res}
R(\widetilde{A}-\l B)\supseteq R(\bpola)\supseteq R(B).
\end{equation}
Now, let $\mu\in\rho(\widetilde{L})$ be arbitrary. From (\ref{eq:res}) follows
\begin{equation*}
(\widetilde{A}-\mu B)x=(\widetilde{A}-\l B)x+(\l-\mu)Bx\in
R(\widetilde{A}-\l B),
\end{equation*}
hence $R(\widetilde{A}-\mu B)\subset
R(\widetilde{A}-\l B)$, which implies $\l\in\rho(\widetilde{L})$.

On the other hand, let
$\l\in\sigma(\widetilde{S})\setminus\sigma_p(\widetilde{S})$.
Then there exists $(y_n)\in\K$ such that
\begin{equation}
\label{eq:ne} (y_n)\notin R(\widetilde{S}-\l).
\end{equation}
We denote $\widetilde{y}=T(y_n)$. Then, the relation
(\ref{eq:ne})
 can be written  as
\[\widetilde{y}\notin R(J\bmin \widetilde{A}-\l\bpola)\]
which implies
\begin{equation}
\label{eq:res2} J\bmin \widetilde{y}\notin
R((\widetilde{A}-\l B)|_{D(\bmin \widetilde{A})}).
\end{equation}
Let us assume that $J\bpola\widetilde{y}\in
R(\widetilde{A}-\l B)$. Then there exists $x\in
D(\widetilde{A})$ such that $\widetilde{A}x-\l Bx=J\bpola
\widetilde{y}$, hence $\widetilde{A}x\in R(\bpola)$, a
contradiction with (\ref{eq:res2}). Hence,
$J\bpola\widetilde{y}\notin R(\widetilde{A}-\l B)$ and
$\l\in\sigma(\widetilde{L})$.
\end{proof}
\begin{proposition}
\label{prop:spec2}
We have
\begin{enumerate}
\item[(i)] $\sigma_c(\widetilde{S})\subset \sigma_c(\widetilde{L})$,
\item[(ii)] $\sigma_{ap}(\widetilde{S})\subset \sigma_{ap}
(\widetilde{L})$.
\end{enumerate}
\end{proposition}
\begin{proof}
 (i) Let $\l\in\sigma_c(\widetilde{S})$. Then
$R(\widetilde{S}-\l)$ is dense in
$\K$. For each $(y_n)\in\K$, and for each $\eps>0$ there
exists $x\in D(\widetilde{S})$ such that
\[
\knorm{(\widetilde{S}-\l)x-(y_n)}<\eps. \] We choose $y\in
R(B)$ arbitrarily and set $\widetilde{y}=J\bmin y$,
$y'=\bmin \widetilde{y}$. For $y'\in\h\subset\K$ there
exists $x\in D(\widetilde{S})$ such that
$\|(\widetilde{S}-\l)x-y'\|_{\K}<\eps\|\bpola\|^{-1}$,
or, equivalently,
\[  \norm{J\bmin\widetilde{A}x-\lambda\bpola x-\widetilde{y}}<\eps\norm{\bpola}^{-1}. \]
Now,
\[  \norm{(\widetilde{A}-\l B)x-y}=\norm{J\bpola(J\bmin
\widetilde{A}x-\l\bpola x-
\widetilde{y})}<\eps. \] Hence, $R(\widetilde{A}-\l B)$ is
dense in $R(B)$, which implies that $R(\widetilde{A}-\l B)$ is dense in
$\h$. Since $\l\notin\sigma_p(\widetilde{L})$,
$\l\in\sigma_c(\widetilde{L})$ follows.

(ii) Let $\in\sigma_{ap}(\widetilde{S})$. Then there
exists a sequence
$(x_n)\in D(\widetilde{S})$, $\knorm{x_n}=1$ such that
$\|(\widetilde{S}-\l) x_n\|\to 0$, or equivalently,
\begin{equation}
\label{eq:ap1} (J\bmin \widetilde{A}-\l\bpola)x_n\to 0
\end{equation}
and
\begin{equation}
\label{eq:ap2} \norm{\bpola x_n}=1.
\end{equation}
Multiplying relation (\ref{eq:ap1}) by $J\bpola$ we obtain
$(\widetilde{A}-\l B)x_n \to 0$. From (\ref{eq:ap2}) we
see that  $\norm{x_n}\ge\norm{\bpola}^{-1}$. Set
$\widehat{x}_n=x_n/\norm{x_n}$. Then
$\norm{\widehat{x}_n}=1$ and $(\widetilde{A}-\l
B)\widehat{x}_n\to 0$, hence
$\l\in\sigma_{ap}(\widetilde{L})$.
\end{proof}
\begin{corollary}
\label{cor:r}
The relation
$\sigma_r(\widetilde{L})\subset\sigma_r(\widetilde{S})$ holds.
\end{corollary}
In general, Proposition \ref{prop:spec2} (and hence also Corollary \ref{cor:r}) cannot be improved, as will be shown
in Example \ref{exmp:Z}.

First we give  sufficient and necessary conditions for the complete correspondence of the spectra
of $\widetilde{L}$ and $\widetilde{S}$.
\begin{lemma}
\label{lem:AA}
The subspace $R(\widetilde{A}-\l B)$ is dense in $\h$ if
and only if $\h\cap R(\widetilde{S}-\l)^{\bot}=\{0\}$.
\end{lemma}
\begin{proof}
$\Longrightarrow$: First we show that the subspace $R(\widetilde{A}-\l B)$ is
dense in $\h$ if and only if $R((\widetilde{A}-\l
B)|_{D(\bmin\widetilde{A})})$ is dense in $\h$.

The case $\l=0$ follows from Lemma
\ref{lem:abc}, so we assume $\l\neq 0$.
First we assume that $R(\widetilde{A}-\l B)$ is
dense in $\h$.  Let $x_0\in\h$, $\eps>0$ be arbitrary.
Then there exists $x\in D(\widetilde{A})$ such that
$\|(\widetilde{A}-\l B)x-x_0\|<\eps/2$. Lemma
\ref{lem:abc} implies that there exists $x'\in
D(\bmin\widetilde{A})$ such that
$\norm{x-x'}<\frac{\eps}{4|\l|\norm{B}}$ and
$\|\widetilde{A}(x-x')\|<\eps/4$. Now,
\[ \norm{(\widetilde{A}-\l B)x'-x_0}\le
\norm{(\widetilde{A}-\l B)(x'-x)+
(\widetilde{A}-\l B)x-x_0}<\eps. \]
The other direction is obvious.

We showed that
$R((\widetilde{A}-\l B)|_{D(\bmin\widetilde{A})})$ is dense in $\h$,
hence the relation
\[ (x,(\widetilde{A}-\lambda B)y)=0, \; \forall y\in D(\bmin \widetilde{A}) \]
implies $x=0$. Now, from
\[ (x,(\widetilde{A}-\lambda B)y)=(\bpola x,\bmin (\widetilde{A}-\lambda B)y)=[x,(\widetilde{S}-\lambda)y]_{\K}  \]
our claim immediately follows.

The other direction follows from the similar reasoning.
\end{proof}
\begin{lemma}
\label{lem:BB}
The subspace $R(\widetilde{S}-\l)$ is dense in $\K$ if and only if \\
$R(\bmin (\widetilde{A}-\l B))$ is dense in $\h$.
\end{lemma}
\begin{proof}
$\Longrightarrow$: Let us assume that $R(\bmin
(\widetilde{A}-\l B))$ is not dense in $\h$. Then there
exists $\widetilde{x}\in\h$, $\widetilde{x}\neq 0$ such
that
\[  (\widetilde{x},\bmin(\widetilde{A}-\l B)y)=0\quad\forall
y\in D(\bmin
\widetilde{A}). \]
Set $(x_n)=T^{-1}J\widetilde{x}$.
Then $((x_n),(\widetilde{S}-\l)y)_{\K}=0$ for all $y\in
D(\widetilde{S})$, which implies the statement. \\
$\Longleftarrow$: Let us assume that $R(\widetilde{S}-\l)$
is not dense in $\K$. Then it exists $(y_n)\in\K$,
$(y_n)\neq 0$, such that
$((y_n),(\widetilde{S}-\l)x)_{\K}=0$ for all $x\in
D(\widetilde{S})$. Set $y= T(y_n)$. Then
\[  (\widetilde{y},(J\bmin \widetilde{A}-\l\bpola)x)=0
\quad \forall x\in
D(\bmin \widetilde{A}), \] hence
\[  (J\widetilde{y},\bmin(\widetilde{A}-\l B)x)=0 \quad
\forall x\in
D(\bmin \widetilde{A}), \] and the statement follows.
\end{proof}
Lemmata \ref{lem:AA} and \ref{lem:BB} imply the following result.
\begin{proposition}
\label{prop:CC}
\begin{enumerate}
\item[(i)] If $\l\in\sigma_r(\widetilde{S})$, then $\l\in\sigma_c
(\widetilde{L})$ if
and only if \\ $\h\cap R(\widetilde{S}-\l)^{\bot}=\{0\}$.
\item[(ii)] If $\l\in\sigma_c(\widetilde{L})$, then $\l\in\sigma_r
(\widetilde{S})$ if
and only if $R(\bmin(\widetilde{A}-\l B))$ is not dense in
$\h$.
\end{enumerate}
\end{proposition}
\begin{example}
\label{exmp:Z}
Let $\h=L_2(0,1)$. We define the operators $A$ and $B$ in $\h$ by:
\begin{gather*}
(Af)(x)=\int_0^1 G(x,\xi)f(\xi)\deri \xi, \\
(Bf)(x)=(2x-1)f(x),
\end{gather*}
where
\[ G(x,\xi)=\begin{cases}
x &, x\le \xi, \\
\xi &, x\ge \xi.
\end{cases} \]
Note that $A$ is the inverse of the Dirichlet Laplacian.
The operator $B$ is bounded and $0\in \sigma_c(B)$, $0\in \sigma_c(A)$ holds. Note that this implies $0\in \sigma_c(L)$.
One can easily prove that $\Pi(L)\cap \Real\ne\emptyset$ holds.

We will show that $R(\bmin A)$ is not dense in $\h$. Indeed, constant functions are not contained in
$R(\bmin A)$. This follows from the fact that the function $x\mapsto \sqrt{|2x-1|}$ is not contained in $R(A)$. Now
Proposition \ref{prop:CC} (ii) implies $0\in\sigma_r(S)$.
\end{example}
\begin{remark}
If $\pi(B)<\infty$ or $\widetilde{A}\ge 0$, then
$\sigma_r(\widetilde{L})=
\sigma_r(\widetilde{S})=\emptyset$, and consequently,
$\sigma_c(\widetilde{L})= \sigma_c(\widetilde{S})$. 

If $\widetilde{S}$ is definitizable, then
$\sigma_r(\widetilde{S})=\emptyset$, hence
$\sigma_r(\widetilde{L})=\emptyset$ and
$\sigma_c(\widetilde{S})= \sigma_c(\widetilde{L})$.
\end{remark}
\begin{corollary}
Let us suppose that $\pi(A)<\infty$. Then each selfadjoint extension
$\widetilde{S}$ of $S$ is definitizable.
\end{corollary}
\begin{proof}
Since deficiency indices of $A$ are finite,
$\pi(\widetilde{A})<\infty$ holds for each selfadjoint extension
$\widetilde{A}$, and Theorem \ref{theorem:kva} implies
$\pi(\widetilde{S})<\infty$. Now, let $\l_0\in\Pi(L)\cap
\Real$ be arbitrary. It is well--known that there
exists a selfadjoint extension $\widetilde{A}_0$ such that
$0\in\rho(\widetilde{A}_0-\l_0 B)$ (see \cite{Kr},
\cite{DM}), which implies $\l_0\in \rho(\widetilde{S}_0)$,
where $\widetilde{S}_0$ is the corresponding operator. Also
$\pi(\widetilde{S}_0)=\pi(\widetilde{A}_0)<\infty$ holds. Now,
from \cite[Proposition 1.1]{CL} follows that each
selfadjoint extension of $S$ has a nonvoid resolvent set,
hence, see \cite{L}, each selfadjoint extension of $S$ is
definitizable.
\end{proof}
\begin{remark}
If $\widetilde{A}$ is qup, the proof of Lemma
3.1 from \cite{BN} implies  $0\in\rho(\widetilde{L})$ for
$\l\neq 0$, $|\l|$ sufficiently small. Hence
$\l\in\rho(\widetilde{S})$ for such $\l$. From this fact
and from $\pi(\widetilde{S})=\pi(\widetilde{A})<\infty$,
it follows that $\widetilde{S}$ is also qup.
\end{remark}
\begin{corollary}
If $\widetilde{L}$ has a compact resolvent in one point
$\l$ (and then, from Proposition \ref{int:prop}, in
all points), and if
\begin{equation}
\label{eq:esch} \widetilde{A}x=\l B x, x\neq 0
\Longrightarrow (Bx,x)\neq 0,
\end{equation}
then $\widetilde{S}$ has a discrete spectrum.
\end{corollary}
\begin{proof}
See \cite{BEL}.
\end{proof}
\begin{remark}
The relation (\ref{eq:esch}) is equivalent with the fact
that $\widetilde{S}$ does not have isotropic eigenvectors, i.e.\ if $[x,x]\ne 0$  for all eigenvectors $x$ of $\widetilde{S}$.
\end{remark}

\section{$B$ unbounded}
\label{sec:unbounded}

In this section we will treat a more general case when
$B$ is possibly unbounded. In this section we assume:
\begin{enumerate}
\item $A$ and $B$ are selfadjoint and $D(A)\subset D(B)$,
\item $0$ is not an eigenvalue of $B$, and
\item $\rho(L)\cap \mathbb{R}\neq \emptyset$.
\end{enumerate}
Using a spectral shift (i.e.\ substituting $A$ by $A-\lambda B$), if necessarily, we can assume
$0\in\rho(A)$. The construction which we use here is essentially
given in \cite{BN}.

Similarly as in Section \ref{sec:bounded}. we introduce a space $\K$ as the
completion of $(D(\bpola),[\cdot,\cdot])$, where
$[\cdot,\cdot]=(J\bpola\cdot,\cdot)$, by the use of Cauchy
sequences. As in Section \ref{sec:bounded}., we set
$J=E(0,\infty)-E(-\infty,0)$, where $E$ is the spectral
function of $B$. By $T$ we denote the operator $T:\K\to
\h$, $T(x_n)=\lim_n \bpola x_n$. It is easy to see that
$T$ is an isomorphism between spaces $\h$ and $\K$. As
was shown in \cite{BN}, $A^{-1}B$ extends to a bounded
symmetric operator $\overline{R}$ in $\K$, and
$\overline{R}$ is injective. Hence, the operator
$S=\overline{R}^{-1}$ is selfadjoint and boundedly
invertible in $\K$. Moreover, we can see that $S$ is an
extension of $B^{-1}A$ which is densely defined in $\K$.
Indeed, since $\overline{R}$ is injective, it holds
\[  (\overline{B^{-1}A})^{-1}=\overline{(B^{-1}A)^{-1}}=
\overline{R}, \]
so $\overline{B^{-1}A}=S$, where the closure is taken in
$\K$. The domain of $S$ is given by
\[   D(S)=\{(x_n)\in\K: J\bmin Ax_n\;
\text{is convergent in}\; \h\}. \]
For $(x_n)\in\K$ we can always assume $x_n\in D(B)$.
Then, if $(y_n)=(S-\l)(x_n)$ we have $T(y_n)=\lim_n
J\bmin(A-\l B)x_n$.

If $A$ is a qup, we can, instead of the assumption 3., make
the assumption $\rho(L)\neq \emptyset$, since in
 \cite[Lemma 3.1]{BN} it was proved that in this case there exists
 $\l_0\in\Real$ such that
$L(\l_0)$ is boundedly invertible. Also, in this case, $S$
is a qup.
\begin{theorem}
\label{theorem:unb} We have
\begin{equation}
\label{eq:unb1} \sigma_p(L)\subset\sigma_p(S)
\end{equation}
 and
\begin{equation}
\label{eq:unb2}
 \sigma(S)\subset\sigma(L).
\end{equation}
\end{theorem}
\begin{proof}
The relation (\ref{eq:unb1}) is obvious.

Let $\l\in\rho(L)$. Then it is easy to see that $\bpola
|A-\l B|^{-1/2}$ and 
$\bpola |A-\overline{\l}B|^{-1/2}$
are everywhere defined bounded operators. From
\[  (|A-\l B|^{-1/2}\bpola)^{\ast}=\bpola
|A-\overline{\l}B|^{-1/2} \]
follows
\[  (|A-\l B|^{-1/2}\bpola)^{\ast\ast}=(\bpola
|A-\overline{\l}B|^{-1/2})^{\ast}, \]
hence
\begin{equation}
\label{eq:unb3}
  |A-\l B|^{-1/2}\bpola \subset (\bpola
  |A-\overline{\l}B|^{-1/2})^{\ast}.
\end{equation}
Since the right hand side in (\ref{eq:unb3}) is a bounded operator, it
follows that also the left hand side in (\ref{eq:unb3}) is a bounded
operator (but, in general, not everywhere defined!). Hence, the operator
\[ \bpola(A-\l B)^{-1}\bpola=\bpola |A-\l B|^{-1/2}\sgn
(A-\l B)|A-\l B|^{-1/2}
\bpola \] is bounded.

The subspace $D(\bpola)\cap
R(\bpola)$ is dense. Indeed, we have
\[
D(\bpola)\cap R(\bpola)=\{x=\bpola y:y\in D(B)\}, \]
 and the right hand side is dense since
$D(B)$ is a core of $\bpola$. Now, let $x\in\h$ be
arbitrary. Then there exists a sequence $(x_n)\in
D(\bpola)\cap R(\bpola)$ such that $x_n\to Jx$. From the
considerations  given above it follows that \\
$\bpola(A-\l B)^{-1}\bpola
x_n$ converges to some $y\in\h$.

Set $z_n=(A-\l
B)^{-1}\bpola x_n$. Then $\bpola z_n\to y$, hence
$(z_n)\in\K$ and $x_n=\bmin(A-\l B)z_n$. Now $\bmin (A-\l
B)z_n\to Jx$, which implies $J\bmin (A-\l B)z_n\to x$, so
$x\in T(R(S-\l))$. This implies $\l\in\rho(S)$.
\end{proof}
The other inclusion $\sigma(L)\subset \sigma(S)$ in
general does not
 hold, since for large $\l$
the operator $L(\l)$ need not be closed.

If we assume $0\in \rho(B)$, then
$\K=D(\bpola)$ and $\sigma_p(S)\subset \sigma_p(L)$, hence
$\sigma_p(S)=\sigma_p(L)$. Also, it can be seen that 
$\sigma_{ap}(L)\subset \sigma_{ap}(S)$. Hence, if $L$ has
only approximative spectrum then $\sigma(L)=\sigma(S)$.

Another kind of construction is made in \cite{LT}. There
it is assumed that $0\in \rho(B)$ and $D(A)=D(\bpola)$.
Then $L(\l)$ is a holomorphic function in Kato sense from
$\K$ to $\h$. For the readers convenience, we give a
definition.
\begin{definition}
 The function $T(\l):\mathcal{X}\to\mathcal{Y}$ is
holomorphic in Kato sence if there exist a bounded
holomorphic functions $U(\l):\mathcal{Z}\to \mathcal{X}$
and $V(\l) :\mathcal{Z}\to \mathcal{Y}$ such that
$T(\l)U(\l)=V(\l)$ and $R(U(\l))=D(T(\l))$ for all
$\l\in\Complex$, where $\mathcal{Z}$ is some Banach space.
An operator function $T(\l)$ is bounded holomorphic if it
takes values in the set of bounded operators and if
$T(\l)x$ is a holomorphic function in $\l$ for all $x$.
\end{definition}
In our case, we can choose $U(\l)=B^{-1}$,
$V(\l)=AB^{-1}-\l$, and $\mathcal{Z}=\h$.

We set $S=B^{-1}A$, and by the use of the similar
methods as in the proofs of Theorem \ref{theorem:spect} and
Proposition \ref{prop:spec2}, it can be seen that
$\sigma(L)=\sigma(S)$, $\sigma_{ap}(L)=\sigma_{ap}(S)$,
$\sigma_c(L)=\sigma_c(S)$ and $\sigma_p(L)=\sigma_p(S)$,
hence we have the complete equivalence of the spectra of $L$ and
$S$.

\section{Variational characterization}
\label{sec:variational}

In this section we apply the constructions of reduced operators from Sections \ref{sec:bounded} and \ref{sec:unbounded} to obtain a Rayleigh--Ritz--like upper bound for the eigenvalues of the reduction operator given in terms of the eigenvalues of the matrix pair $(A_V,B_V)$, where $A_V$, $B_V$ are orthogonal projections of $A$, $B$ on a finite dimensional subspace. 
Hence the bound is given in terms of the original data $(A,B)$, which can be important in applications.

To obtain a variational characterization for the reduction operator, we need a result from \cite{BN94}.

Let $Q$ be a qup operator such that all finite critical
points are regular, and such that there are no critical
points embedded in the positive continuous spectrum. For such a
$Q$, in \cite{BN94} a variational characterization of
eigenvalues is obtained. More precisely, we define
\begin{equation}
\label{eq:var}
  \l_j^{\pm}=\sup\left\{\inf\left\{\frac{[Qx,x]}{[x,x]}: x\in M
  \cap\mathcal{C}^{\pm}\cap
D(Q)\right\}:M\in\mathcal{M}_j\right\},
\end{equation}
where $\mathcal{C}^{\pm}=\{x: [x,x]\gtrless 0\}$, and
$\mathcal{M}_j$ denotes the set of all subspaces of
codimension $j-1$. Then, if by $d^{\pm}$ we denote the
positive (negative) spectral shift of $Q$ (for more
details see \cite{BN94}), a number $\l_{j+d^{\pm}}^{\pm}$
is either an eigenvalue of $Q$, or a point on the boundary
of the essential spectrum of $Q$, for all $j\in\mathbb{N}$.

If the operator $\widetilde{S}$ from Section \ref{sec:bounded}, or the operator $S$ from Section \ref{sec:unbounded} satisfy conditions
stated above, we can variationally characterize the eigenvalues of $ \widetilde{S}$ or $S$, respectively.
If $B$ is bounded, then we have the following  variational characterization
\begin{multline}
\label{eq:var2}
  \l_j^{\pm}(\widetilde{S})=\l_j^{\pm}(\widetilde{L})\\
  =\sup \left\{\inf\left\{\frac{(\widetilde{A}x,x)}{(Bx,x)}: x\in M
  \cap\mathcal{C}^{\pm}\cap
D(\bmin \widetilde{A})\right\}:M\in\mathcal{M}_j\right\}.
\end{multline}
In the case when $B$ is unbounded, we have
\[  \frac{[S(x_n),(x_n)]_{\K}}{[(x_n),(x_n)]_{\K}}=
\lim_n \frac{(Ax_n,x_n)}{(Bx_n,x_n)}, \] where
$(x_n)\in D(S)$. Since $D(\bmin A)$ is a core of the operator $S$,
from \cite[Theorem 4.5.3.]{Davies} follows that in
this case we have
\begin{multline}
\label{eq:var2-2}
  \l_j^{\pm}(S)=\l_j^{\pm}(L)\\
  =\sup \left\{\inf\left\{\frac{(Ax,x)}{(Bx,x)}: x\in M
  \cap\mathcal{C}^{\pm}\cap
D(\bmin A)\right\}:M\in\mathcal{M}_j\right\}.
\end{multline}
Now, suppose that the operator $ \widetilde{S}$ from Section \ref{sec:bounded} satisfies conditions stated above.
Let $V\subset D(\bmin \widetilde{A})$ be a
finite dimensional ortho--complemented subspace in $\K$. Let
$\dim V=N$. Denote by $P$ a orthogonal projector in $\K$
onto $V$. Set $\widetilde{A}_V=P\widetilde{A}P$,
$B_V=PBP$. Let $\bot$ and $[\bot]$ denote the  orthogonal complement in
$\h$ and $\K$, respectively. The set of
ortho--complemented subspaces of $\K$ we denote by
$\mathcal{O}$. Let $F$ be a subspace of $V$, with $\dim F=k$, and let
$F=\spa\{x_1,\ldots, x_k\}$. Set
$F'=\spa\{B_V^{-1}x_1,\ldots,B_V^{-1}x_k\}$. Then $\dim
F'=k$ and
\[   x [\bot] F' \Longleftrightarrow x\bot F, \quad
\forall x\in V. \]
On the other hand, if $x [\bot] F$, then we define
$F''=\spa\{B_Vx_1,\ldots, B_V x_k\}$ and then we have
\[ x [\bot] F \Longleftrightarrow x\bot F'',\quad\forall x\in V. \]
For $m\le N$, set
\begin{equation}
\label{eq:mat}
  \mu_m=\sup_{\substack{
  F\subset V \\
  \dim F=m-1}}
  \inf_{\substack{
  x\in V\\
  x\bot F \\
  x\in\mathcal{C}^+}}
  \frac{(\widetilde{A}_Vx,x)}{(B_Vx,x)}.
  \end{equation}
Then, from the considerations given above, we have
\begin{equation*}
\begin{split}
\mu_m &=\sup_{\substack{
  F\subset V \\
  \dim F=m-1}}
  \inf_{\substack{
  x\in V\\
  x[\bot] F \\
  x\in\mathcal{C}^+}}
  \frac{(\widetilde{A}x,x)}{(Bx,x)} \\
  &=\sup_{\substack{
  y_1,\ldots,y_{m-1}\in V \\
  \spa\{y_1,\ldots,y_{m-1}\}\in\mathcal{O}}}
  \inf_{\substack{
  x\in V\\
  x\in\spa\{y_1,\ldots,y_{m-1}\}^{[\bot]} \\
  x\in\mathcal{C}^+}}
  \frac{(\widetilde{A}x,x)}{(Bx,x)} \\
&=\sup_{\substack{
  y_1,\ldots,y_{m-1}\in \K \\
  \spa\{y_1,\ldots,y_{m-1}\}\in\mathcal{O}}}
  \inf_{\substack{
  x\in V\\
  x\in\spa\{Py_1,\ldots,Py_{m-1}\}^{[\bot]} \\
  x\in\mathcal{C}^+}}
  \frac{(\widetilde{A}x,x)}{(Bx,x)} \\
  &=\sup_{\substack{
  y_1,\ldots,y_{m-1}\in \K \\
  \spa\{y_1,\ldots,y_{m-1}\}\in\mathcal{O}}}
  \inf_{\substack{
  x\in V\\
  x\in\spa\{y_1,\ldots,y_{m-1}\}^{[\bot]} \\
  x\in\mathcal{C}^+}}
  \frac{(\widetilde{A}x,x)}{(Bx,x)} \ge
\end{split}
\end{equation*}
\begin{equation*}
\begin{split}
  &\ge \sup_{\substack{
  y_1,\ldots,y_{m-1}\in \K \\
  \spa\{y_1,\ldots,y_{m-1}\}\in\mathcal{O}}}
  \inf_{\substack{
  x\in D(\bmin \widetilde{A})\\
  x\in\spa\{y_1,\ldots,y_{m-1}\}^{[\bot]} \\
  x\in\mathcal{C}^+}}
  \frac{(\widetilde{A}x,x)}{(Bx,x)} \\
  &=\sup_{\substack{
  M\in \mathcal{M}_j}}
  \inf_{\substack{
   x\in D(\bmin \widetilde{A})\\
  x\in M\cap\mathcal{C}^+}}
  \frac{(\widetilde{A}x,x)}{(Bx,x)}.
\end{split}
\end{equation*}
The analogous relation can be obtained for the eigenvalues of
the negative type. From the finite dimensional principle given in
\cite{BNY}, $\mu_m$ is an eigenvalue of the matrix pair
$(\widetilde{A}_V, B_V)$. Hence we obtained a
Rayleigh--Ritz--like upper bound for eigenvalues of
$\widetilde{S}$:
\[  \l_j^{\pm}(\widetilde{A}_V,B_V)\ge \l_j^{\pm}(\widetilde{S}). \]
The analogous Rayleigh--Ritz upper bound holds also for the operator $S$ from Section \ref{sec:unbounded},
due to the formula (\ref{eq:var2-2}).

\section{Zero is an eigenvalue of $B$}
\label{sec:zero}

If $0$ is an eigenvalue of $B$, then in general a reduced operator cannot be constructed. But we can construct an reduced linear relation, a generalization of the notion of linear operator, and we can recover the results from previous sections. 
Under some additional assumptions we can also construct a proper reduction operator.

In this section we assume:
\begin{enumerate}
\item $A$ and $B$ are selfadjoint operators,
\item $B$ is $A$--bounded, $0\in\sigma_p(B)$,
\item $\rho(A,B)\cap \Real\neq \emptyset$.
\end{enumerate}
The last relation implies that we can assume $0\in\rho(A)$.

Then $B^{-1}A$ does not exist as a linear operator, but
it can be introduced in terms of the linear relations, i.e.\ as
a subspace of $\h\times\h$ (for the basic definitions see
\cite{S}, \cite{DS}). We define
\[  S_0=B^{-1}A=\{\{x,y\}\in D(A)\times D(B): Ax=By\}. \]
By $S$ we denote the closure of $S_0$ as a subspace in
$\h\times\h$.

Let $\K$ denote the factor space $D(B)/_{\Ker B}$, and let
$\widetilde{\K}$ denote the completion of
$(\K,(B\cdot,\cdot))$. Evidently, $\widetilde{\K}$ is a
Krein space, and it consists of all Cauchy sequences in
$D(\bpola)$, where we identify such two sequences $(x_n)$
and $(y_n)$ if $\bpola (x_n-y_n)\to 0$. This equivalence relation
 we denote by the $\sim$ symbol. Let $[x]$ denote a
class in $\widetilde{\K}$ represented by $x$.

Let $R(\l,S)$ be a resolvent of $S$. Then the family of
linear operators $R(\l,S)$ in $\h$ induces a family
$\widetilde{R}(\l,S)$ of linear operators in the space $\K$,
since $R(\l,S)\Ker B\subset\Ker B$ (see \cite{LT}). We
define the linear relation $\widetilde{S}$ in
$\widetilde{\K}$ by $\widetilde{S}=\widetilde{R}(0,S)^{-1}$.

We say that $\l$ is a \textsl{discrete eigenvalue} of the
operator $T$ if $\l$ is an isolated eigenvalue of finite
multiplicity of the operator $T$.

First we treat the case when $B$ is a bounded
operator and zero is a discrete eigenvalue. Then it is easy
to see that we can substitute Cauchy sequences with
elements of $\h$, and that $x\sim y$ if and only if
$x-y\in\Ker B$. The inner product on $\widetilde{K}$ is
given by $[\cdot,\cdot]=(B\cdot,\cdot)$. It is easy to
see that
\[  \widetilde{S}=\{\{[x],[y]\}:x-A^{-1}By\in\Ker B\}, \]
and that $\widetilde{S}$ is a selfadjoint
relation.
\begin{theorem}
The spectra of $\widetilde{S}$ and $L$ coincide. Moreover,
their point spectra coincide.
\end{theorem}
\begin{proof}
The relation $\sigma_p(L)\subset \sigma_p(\widetilde{S})$
is obvious. On the other hand, let \\
$\l\in\sigma_p(\widetilde{S})$. Then there exists $x\in\h$
such that $x-\l A^{-1}Bx\in\Ker B$, hence $B x-\l BA^{-1}B
x=0$. There exists $y\neq 0$ such that $B x= A y$, hence
\[ 0=(I-\l B A^{-1})Bx=(I-\l B A^{-1})Ay=Ay-\l B y, \]
which implies $\l\in\sigma_p(L)$. It is easy to see that
\[   (S-\l)^{-1}=(A-\l B)^{-1} B\quad\forall\l\in
\Complex\setminus\sigma_p(L)=
\Complex\setminus\sigma_p(\widetilde{S}). \] From this
relation it is easy to obtain
$\rho(L)=\rho(\widetilde{S})$.
\end{proof}

If $B$ is bounded, but zero is not a discrete eigenvalue of
$B$, it can be seen that
\[ \widetilde{S}=\{\{[x],[(y_n)]\}:  x- A^{-1}J\bpola\lim_n
\bpola y_n \in \Ker B \}, \] 
where $J$ is defined as in section
\ref{sec:bounded}, is a selfadjoint relation. Using similar techniques as
in Section \ref{sec:bounded}, it can be shown that
$\sigma(\widetilde{S})=\sigma(L)$ and
$\sigma_p(\widetilde{S})=\sigma_p(L)$.

This kind of procedure can be implemented also in the case
when $B$ is not bounded, using the techniques from Section
\ref{sec:unbounded}. For instance, if all assumptions given in
the beginning of the Section \ref{sec:unbounded} hold, except from
$0\notin\sigma_p(B)$, we can define $\widetilde{S}$ by
\[  \widetilde{S}=\{\{[(x_n)],[(y_n)]\}:(x_n)\sim (A^{-1}By_n)\}. \]
Reasoning analogously as in the proof of Theorem \ref{theorem:unb}, it can be shown that $\sigma_p(L)\subset \sigma_p(\widetilde{S})$ and $\sigma(\widetilde{S})\subset \sigma(L)$.

When is it possible to reduce the relation $\widetilde{S}$ to
an operator, without the loss of the information about spectra?
To answer this question we introduce the notion of the
\textsl{multi--valued part}
 $T_{\infty}=\{\{0,g\}\in T\}$ of the relation $T$ (see \cite{S}).
We also set $T(0)=\{g\in\K:\{0,g\}\in T\}$. It is easy to
see that $T_{\infty}$ is ortho--complemented if and only if
$T(0)$ is ortho--complemented. We have
\begin{theorem}[\cite{S}]
\label{theorem:S} Let $T$ be a closed symmetric relation in a
Krein space $\K$ such that $T_{\infty}$ is
ortho--complemented. Then $T_s=T\cap (T_{\infty})^{\bot}$ is
an operator with $D(T_s)=D(T)$ and $R(T_s)\subset
T(0)^{\bot}$. We also have
\begin{enumerate}
\item[(i)] $\sigma_p(T)=\sigma_p(T_s)$,
\item[(ii)] $\sigma(T)=\sigma(T_s)$, and
\item[(iii)] $T$ is a selfadjoint relation if and only if $T_s$ is a
selfadjoint
operator in $T(0)^{\bot}$.
\end{enumerate}
\end{theorem}
Hence, we can reduce $\widetilde{S}$ to an operator if
the set
\[ \widetilde{S}(0)=\{[(y_n)]: \lim_n \bpola A^{-1}B y_n \in
\Ker B\} \]
is ortho--complemented.

If, for instance, $R(B)\cap A\Ker B=\{0\}$, or if $\Ker B$
is $A$--invariant, $\widetilde{S}(0)$ is
ortho--complemented.

If $d_0=\dim\Ker B<\infty$ and if $\widetilde{S}(0)$
is ortho--complemented, we can also obtain
Rayleigh--Ritz--like upper bound for the eigenvalues of
$\widetilde{S}(0)$ in terms of the non-zero eigenvalues of
the matrix pair $(A_V,B_V)$, where $V\subset D(A)$, in the case when the
dimension of $V$ is greater then the sum of $d_0$ and the
spectral shift of $\widetilde{S}_s=\widetilde{S}\cap
(\widetilde{S}(0))^{\bot}$.
We omit the details.

\medskip

\noindent\textbf{Acknowledgement}: The author is very grateful to anonymous referees for extremely useful comments and suggestions.

\bibliography{bmina}
\bibliographystyle{plain}
\end{document}